\magnification=1080
\vsize=192mm
\hsize=135mm
\voffset=5mm
\null


\def\Pb{1}
\def\Ha{2}
\def\Hb{3}
\def\Hc{4}
\def\Hd{5}
\def\He{6}
\def\Hf{7}
\def\Hg{8}
\def\Hh{9}
\def\La{10}
\def\Lb{11}
\def\Lc{12}
\def\Da{13}
\def\Daa{14}
\def\Ca{15}
\def\Cb{16}
\def\Cbb{17}
\def\Db{18}
\def\Cc{19}
\def\Ld{20}

\def\BDS{1}
\def\CDF{2}
\def\CGV{3}
\def\CV{4}
\def\DF{5}
\def\DFPV{6}
\def\FLaa{7}
\def\FPT{8}
\def\FT{9}
\def\GV{10}
\def\KaA{11}
\def\KaB{12}
\def\MHM{13}
\def\Mo{14}
\def\Pi{15}
\def\PSY{16}
\def\PSZ{17}
\def\QSb{18}

\null

\newtoks\hautpagegauche
\newtoks\hautpagedroite
\newtoks\paragraphecourant
\newtoks\chapitrecourant
\hautpagegauche={}
\hautpagedroite={}
\headline={\ifodd\pageno\the\hautpagedroite\else\the\hautpagegauche\fi}

\font\TenEns=msbm10
\font\SevenEns=msbm7
\font\FiveEns=msbm5
\newfam\Ensfam
\def\Ens{\fam\Ensfam\TenEns}
\textfont\Ensfam=\TenEns
\scriptfont\Ensfam=\SevenEns
\scriptscriptfont\Ensfam=\FiveEns
\def\R{{\Ens R}}

\def\Rn{{\R}^n}


\font\itsmall=cmsl9

\font\eightrm=cmr9
\font\sixrm=cmr6
\font\fiverm=cmr5

\font\eighti=cmmi9
\font\sixi=cmmi6
\font\fivei=cmmi5

\font\eightsy=cmsy9
\font\sixsy=cmsy6
\font\fivesy=cmsy5

\font\eightit=cmti9
\font\eightsl=cmsl9
\font\eighttt=cmtt9

\def\eightpoint{\def\rm{\fam0\eightrm}
\textfont0=\eightrm
\scriptfont0=\sixrm
\scriptscriptfont0=\fiverm

\textfont1=\eighti
\scriptfont1=\sixi
\scriptscriptfont1=\fivei

\textfont2=\eightsy
\scriptfont2=\sixsy
\scriptscriptfont2=\fivesy

\textfont3=\tenex
\scriptfont3=\tenex
\scriptscriptfont3=\tenex

\textfont\itfam=\eightit \def\it{\fam\itfam\eightit}
\textfont\slfam=\eightsl \def\it{\fam\slfam\eightsl}
\textfont\ttfam=\eighttt \def\it{\fam\ttfam\eighttt}
}


\font\pc=cmcsc9
\font\itsmall=cmsl9

\def \trait (#1) (#2) (#3){\vrule width #1pt height #2pt depth #3pt}
\def \fin{\hfill
       \trait (0.1) (5) (0)
       \trait (5) (0.1) (0)
       \kern-5pt
       \trait (5) (5) (-4.9)
       \trait (0.1) (5) (0)
\medskip}

\input color

\paragraphecourant={\rmt }
\chapitrecourant={\rmt Ph. SOUPLET}
\footline={\ifnum\folio=1 \hfill\folio\hfill\fi}
\hautpagegauche={\tenrm\folio\hfill\the\chapitrecourant\hfill}
\hautpagedroite={\ifnum\folio=1 \hfill\else\hfill\the\paragraphecourant\hfill\tenrm\folio\fi}

\font\rmb=cmbx8 scaled 1125 \rm

\centerline{\rmb GLOBAL EXISTENCE FOR REACTION-DIFFUSION SYSTEMS}
\vskip 0.5mm
\centerline{\rmb WITH DISSIPATION OF MASS AND QUADRATIC GROWTH}

\vskip 5mm
\centerline{\pc Philippe SOUPLET}
\vskip 3mm
\centerline{\itsmall Universit\'e Paris 13, Sorbonne Paris Cit\'e, CNRS UMR 7539}
\centerline{\itsmall Laboratoire Analyse G\'eom\'etrie et Applications}
\centerline{\itsmall 93430 Villetaneuse, France. Email: souplet@math.univ-paris13.fr}

\vskip 4mm

\baselineskip=12pt
\font\rmt=cmr9

\setbox1=\vbox{
\hsize=120mm
{\baselineskip=11pt \parindent=3mm \eightpoint \rmt
{\pc Abstract:}\ We consider the Neumann and Cauchy problems for positivity 
preserving reaction-diffusion systems of $m$ equations enjoying the mass and entropy dissipation properties.
We show global classical existence in any space dimension, under the assumption that the nonlinearities have at most
quadratic growth. This extends previously known results which, in dimensions $n\ge 3$, required mass conservation
and were restricted to the Cauchy problem. Our proof is also simpler.
\vskip 0.2cm

{\pc Keywords:}\ reaction-diffusion systems, mass dissipation, entropy, global existence

}}

\hskip 2mm \hbox{\box1}

\bigskip

{\bf 1. Introduction and main result.}
\medskip

We consider the reaction-diffusion system
$$
\left\{\eqalign{
\partial_tu_i-d_i\Delta u_i &=f_i(u),\quad x\in\Omega,\ 0<t<T\quad (1\le i\le m),\cr
\partial_\nu u_i &=0,\quad x\in\partial\Omega,\ 0<t<T,\cr
u_i(x,0) &=u_{i,0}(x),\quad x\in\Omega,}\right.
\leqno(\Pb)$$
where $\Omega\subset \Rn$ is either a smoothly bounded domain, or $\Omega=\Rn$, with $n\ge 1$.
Here $m\ge 2$, $d_i>0$, $u=(u_1,\dots,u_m)$, and the functions $f_i: \R^m\to \R$ are of class $C^1$ for  $i=1,\dots,m$.
As for the initial data, we assume
$$\hbox{ $u_0\in (L^\infty(\Omega))^m$, \ with $u_0\ge 0$ in $\overline\Omega$ and $u_{0,i}\not\equiv 0$ for  $i=1,\dots,m$} 
\leqno(\Ha)$$
(in this paper vector inequalities such as $u\ge 0$ or $u>0$ will be understood component-wise).
Then problem (\Pb) admits a unique, maximal solution,
classical for $t>0$. Its existence time will be denoted by~$T_{max}$.

We assume the following structure conditions:
$$\hbox{(Preservation of positivity)}\qquad\qquad u\ge 0,\ u_i=0\Longrightarrow f_i(u)\ge 0,
\leqno(\Hb)$$
$$\hbox{(Dissipation of mass)}\qquad\qquad\sum_{i=1}^m f_i(u)\le 0,\quad u>0.
\leqno(\Hc)$$
Assumption (\Hb) implies that $u>0$ in $\overline\Omega\times (0,T_{max})$ by the maximum principle. 
Assumption (\Hc) guarantees that the total mass $\sum_i\int_\Omega u_i(x,t)\, dx$ is nonincreasing in time
(in case $\Omega$ is bounded), 
and such systems are frequently encountered in models of chemical reactions.
In the equidiffusive case $d_i\equiv d$, assumptions (\Hb)(\Hc) ensure global existence and boundedness of the solution,
as an immediate consequence of the maximum principle applied to $\sum_i u_i$.

It is well known that global existence is no longer trivial when the $d_i$ are not equal
(a case which is indeed relevant in chemical reactions)
and there is an abundant mathematical literature on this question (see, e.g., 
[\QSb, Section 33] and [\Pi] for surveys, and see also further references in Remark 1 below).
Various sufficient conditions on the nonlinearities $f_i$ for global classical existence have been found,
as well as examples of finite time blowup for certain systems.
The case of systems with at most quadratic growth 
$$|f_i(u)|\le M(1+|u|^2),\quad u\ge 0,\ 1\le i\le m,
\leqno(\Hd)$$
has received particular interest,
especially in view of the special case 
$$m=4,\quad f_i(u)=(-1)^{i}(u_1u_3-u_2u_4),
\leqno(\He)$$
which corresponds to the reversible binary reaction
$$A_1+A_3 \buildrel {\displaystyle\longleftarrow}\over {\longrightarrow} A_2+A_4.$$
It was proved by Kanel [\KaB] that if $\Omega=\Rn$, (\Ha), (\Hd) are satisfied and assumption (\Hb)
is replaced with the stronger mass conservation assumption
$$\sum_{i=1}^m f_i(u)=0,\quad u\ge 0
\leqno(\Hf)$$
(which guarantees that the total mass $\sum_i\int_\Omega u_i(x,t)\, dx$ is conserved in time),
then global existence is true in any space dimension. This in particular covers the case~(\He).
See also [\KaA, \CV, \GV, \CDF, \CGV] for related results and alternative approaches.
In particular the case (\He) in a bounded domain in dimension $2$ is covered in [\GV].

The goal of the present paper is to extend the global existence result 
for system~(\Pb) with at most quadratic growth, under the weaker mass dissipation condition (\Hb),
as well as to cover the case of bounded domains with Neumann conditions.
As a counterpart, we will need an additional structure assumption, the so-called dissipation of entropy property:
$$\sum_{i=1}^m f_i(u)\log u_i\le 0,\quad u>0.
\leqno(\Hg)$$
Such a condition is satisfied by many systems corresponding to reversible reactions,
and this is for instance true in the case of (\He).
As for the at most quadratic growth condition, we will assume it under the following slightly stronger form
$$|\nabla f_i(u)|\le M(1+|u|),\quad u\ge 0,\ 1\le i\le m,
\leqno(\Hh)$$
for some constant $M>0$. 
Our main result is the following.

\proclaim Theorem.
Let $m\ge 2$, $n\ge 1$, let $\Omega\subset \Rn$ be either a smoothly bounded domain, or $\Omega=\Rn$,
 and let $u_0$ satisfy (\Ha).
Assume that the nonlinearities $f_i\in C^1(\R,\R)$ satisfy properties (\Hb), (\Hc), (\Hg), (\Hh).
Then problem (\Pb) has a global classical solution, i.e. $T_{max}=\infty$.

Our proof is based on suitable modifications of Kanel's methods in [\KaA, \KaB].
Namely, we combine the argument from [\KaA], based on interpolation inequalities and suitable auxiliary problems,
 with the entropy structure provided by (\Hg).
This enables one to reach nonlinearities with quadratic growth
under a mere mass dissipation structure, without requiring mass conservation.
Actually, the mass conservation structure is crucially used in [\KaB] and [\CGV], via delicate H\"older estimates
for suitable parabolic equations with bounded coefficients (based on De Giorgi type iteration).
We can here take advantage of the entropy structure to avoid such difficulties,
resulting in a much simpler argument.
We are also able to extend this approach to the Neumann problem.

\bigskip

{\bf Remark 1.}
(a) Theorem~1 remains valid, with simple proof changes, if one considers 
Dirichlet instead of Neumann boundary conditions in (\Pb).

\smallskip

(b) Assumption (\Hh) in Theorem~1 can be replaced by the slightly weaker condition:
$$|f_i(u)|\le M(1+|u|^2)\quad\hbox{and}\quad
\textstyle{\partial f_i\over \partial u_j}(u)\ge -M(1+|u|),\quad u\ge 0,\ 1\le i,j\le m.
\leqno(\Hh')$$
Also, assumption (\Hg) can be replaced with
$$\sum_{i=1}^m f_i(u)(1+\log u_i)\le C\sum_{i=1}^m u_i \log (1+u_i),\quad\hbox{ for all $u$ with $u_i\ge 1$, $1\le i\le m$}.
$$

(c) In the case $\Omega$ is bounded, conditions (\Hb) and (\Hc) guarantee that $u$ remains bounded in $L^1$.
But it is an open problem whether $u$ is globally bounded in $L^\infty$
under the assumptions of Theorem~1 with $n\ge 3$, even for the case of system (\He).
For positive results when $n\le 2$, see [\PSY] and the references therein.
\smallskip

(d) The existence of a global {\it weak} solution was shown in [\DFPV, \Pi]
under the mere assumptions (\Hb)--(\Hd) (without conservation of mass or dissipation of entropy).
\smallskip

(e) A related, active topic is the study, by means of entropy methods, of the stabilization as $t\to\infty$
of (classical or weak) global solutions of systems with dissipation of mass.
For this, we refer to, e.g., [\DF, \BDS, \MHM, \FLaa, \FPT, \FT, \PSZ, \PSY]. 

\bigskip
{\bf 2. Proof of Theorem 1.}
\medskip

We shall use the following interpolation lemma. It was proved in [\KaA] in the case $\Omega=\Rn$
and we here extend it to the case of bounded domains with Neumann boundary conditions.
The proof will be given in Section~3.

For $T>0$, we denote $Q_T=\Omega\times (0,T)$. For $k=1,2$, we set $E_k=\{\psi\in BC^k(\overline\Omega) : \hbox{$\partial_\nu \psi=0$ on $\partial\Omega$}\}$
(where the boundary conditions is omitted if $\Omega=\Rn$, and where $\psi$ may be real- or vector-valued).
We also denote $\|U\|_{k,T}=\sup_{t\in (0,T)}\|U(t)\|_{C^k(\overline\Omega)}$ for $k\ge 0$ integer and,
if $u$ is vector-valued, $\|u\|_{k,T}=\max_{1\le i\le m}\|u_i\|_{k,T}$.

\proclaim Lemma~2.
Let $T>0$, $U_0\in E_1$, $g\in BC(\overline Q_T;\R)$ and let $U:\overline Q_T\to\R$ be a classical solution~of
$$
\left\{\eqalign{
U_t-d\Delta U &=g,\quad x\in\Omega,\ 0<t<T,\cr
U_\nu &=0,\quad x\in\partial\Omega,\ 0<t<T,\cr
U(x,0) &=U_0(x),\quad x\in\Omega.}\right.
$$
(i) Then we have
$$\|U\|_{1,T} \le C(\Omega,d,T)\Bigl[\|U_0\|_{C^1}+\|U\|_{0,T}^{1/2}\|g\|_{0,T}^{1/2}\Bigr].
\leqno(\La)$$
Here and below, the constants $C(\Omega,d,T)>0$ remain bounded for $T>0$ bounded.
\smallskip
(ii) Assume in addition that $U_0\in E_2$ and $g\in C([0,T];E_1)$.
Then we have
$$\|U\|_{2,T} \le C(\Omega,d,T)\Bigl[\|U_0\|_{C^2}+\|U\|_{1,T}^{1/2}\|g\|_{1,T}^{1/2}\Bigr]
\leqno(\Lb)$$
and
$$\|U\|_{2,T} \le C(\Omega,d,T)\Bigl[\|U_0\|_{C^2}+\|g\|_{1,T}^{1/2}
\Bigl[\|U_0\|_{C^1}+\|U\|_{0,T}^{1/2}\|g\|_{0,T}^{1/2}\Bigr]^{1/2}\Bigr].
\leqno(\Lc)$$

\medskip

{\bf Proof of Theorem 1.} By a time shift we may assume without loss of generality that 
$u_0\in BC^2(\overline\Omega)$, with $\partial_\nu u_0=0$ on $\partial\Omega$ if $\Omega$ is bounded.
Also we shall write $\sum_i$ for $\sum_{i=1}^m$.
\medskip

{\it Step 1.\ Passage to entropy variables.}
We set $L_i=\partial_t-d_i\Delta$ and define the new unknowns 
$$v_i:=(1+u_i)\log(1+u_i)>0,\qquad w_i:=v_i e^{-Kt}.$$
We claim that for suitable constant $K>0$, the functions $w_i$ satisfy
$$\sum_i L_iw_i\le 0.
\leqno(\Da)$$

To this end we compute
$$\partial_t v_i=(1+\log(1+u_i))\partial_tu_i,\quad
\nabla v_i=(1+\log(1+u_i))\nabla u_i$$
and
$$
\Delta v_i=(1+\log(1+u_i))\Delta u_i+(1+u_i)^{-1}|\nabla u_i|^2,$$
hence
$$L_iv_i\le (1+\log(1+u_i))L_iu_i=(1+\log(1+u_i))f_i(u).$$
Set $e=(1,\dots,1)$ and denote by $|\cdot|_\infty$ the max norm on $\R^m$.
It follows from (\Hc), (\Hg), (\Hh) (or (\Hh')) and the mean value theorem that
$$\eqalign{
\sum_i L_iv_i
&\le \sum_i \log(1+u_i)f_i(u)\cr
&= \sum_i \log(1+u_i)(f_i(u)-f_i(e+u))+\sum_i \log(1+u_i)f_i(e+u) \cr
&\le \sum_i \log(1+u_i)(f_i(u)-f_i(e+u))
\le m^{3/2}M (1+|u|_\infty)\log(1+|u|_\infty) \cr
&\le m^{3/2}M\sum_i (1+u_i)\log(1+u_i)=m^{3/2}M\sum_i v_i.
}$$
We deduce (\Da) with $K=m^{3/2}M$.
\smallskip

{\it Step 2.\ Linear auxiliary problem.}
Pick any finite $T<T_{max}$. Following [\KaA, \KaB]
(in slightly modified form), we fix $d=1+\max_i d_i$, set $L=\partial_t-d\Delta$ 
and, for each $1\le i\le m$, we introduce the (classical) solution $z_i\ge 0$ 
of the auxiliary problem
$$
\left\{\eqalign{
Lz_i &=w_i,\quad x\in\Omega,\ 0<t<T,\cr
\partial_\nu z_i &=0,\quad x\in\partial\Omega,\ 0<t<T,\cr
z_i(x,0) &=0,\quad x\in\Omega.}\right.
\leqno(\Daa)$$
We claim that there exists a constant $C_1>0$ 
independent of $T$ such that
$$w_i\le C_1-\sum_i(d-d_i)\Delta z_i \quad\hbox{in $Q_T$ for $1\le i\le m$}
\leqno(\Ca)$$
and
$$z_i\le dC_1T\quad\hbox{in $Q_T$ for $1\le i\le m$.}
\leqno(\Cb)$$

To show (\Ca) and (\Cb), we set 
$$\phi=\sum_i L_iz_i.$$
Using (\Da), we first notice that $\phi$ satisfies
$$L\phi=\sum_i L_i(Lz_i)=\sum_i L_iw_i\le 0.
\leqno(\Cbb)$$
Next, in the case $\Omega$ bounded, 
observing that $z_i$ 
 is sufficiently smooth up to the boundary and that 
$\partial_\nu\partial_t  z_i=\partial_t\partial_\nu  z_i =0$ on $\partial\Omega$, we get
$$\partial_\nu \Delta z_i=-d^{-1}\partial_\nu Lz_i=-d^{-1}\partial_\nu w_i=0 \quad\hbox{ on $\partial\Omega$.}$$
Therefore, 
$\partial_\nu(L_iz_i)=-d_i\partial_\nu \Delta z_i=0$, hence $\partial_\nu \phi=0$, on $\partial\Omega$.
In both cases $\Omega$ bounded and $\Omega=\Rn$, it thus follows from (\Cbb) and the maximum principle that 
$$\phi\le C_1\quad\hbox{in $Q_T$,}
\leqno(\Db)$$
with $C_1$ independent of $T$.
Then eliminating $\partial_tz_i$ between $L$ and $L_i$, by writing
$$\sum_i(d-d_i)\Delta z_i=\sum_i(L_iz_i-Lz_i)=\phi-\sum_i w_i,$$
we see that (\Db) and $w_i\ge 0$ guarantee (\Ca).

On the other hand, we may eliminate $\Delta z_i$ by writing
$$\sum_i(d-d_i)\partial_tz_i=\sum_i(dL_iz_i-d_iLz_i)=d\phi-\sum_i d_iw_i\le dC_1.$$
Integrating in time and using $z_{i,0}=0$ and $d=1+\max_i d_i$, we get
$$\sum_iz_i\le \sum_i(d-d_i)z_i\le \sum_i(d-d_i)z_{i,0}+dC_1T\le dC_1T,$$
hence (\Cb), owing to $z_i\ge 0$.
\smallskip

{\it Step 3.\ Interpolation and feedback argument.}
We shall now use inequalities (\Ca), (\Cb), along with a feedback argument, to bound $w_i$ (hence $u_i$).
To this end, we shall suitably estimate the diffusion terms $\Delta z_i$ by means of the interpolation Lemma~2.
In this~Step~3, $C(T)$ will denote a generic positive constant (possibly depending on the solution $u$),
which remains bounded for $T>0$ bounded.

By  (\Ca), (\Daa), (\Lc) in Lemma~2 and (\Cb), we have
$$\eqalign{
\|w_i\|_{0,T}
&\le C\bigl[1+\|z_i\|_{2,T} \bigr]
\le C(T)\Bigl[1+\|w_i\|_{1,T}^{1/2}\,\|z_i\|_{0,T}^{1/4}\,\|w_i\|_{0,T}^{1/4}\Bigr] \cr
\noalign{\vskip 0.2mm}
&\le C(T)\Bigl[1+\|w_i\|_{1,T}^{1/2}\,\|w_i\|_{0,T}^{1/4}\Bigr],}$$
hence
$$\|w_i\|_{0,T}\le C(T)\Bigl[1+\|w_i\|_{1,T}^{2/3}\Bigr].
\leqno(\Cc)$$
On the other hand, since $|f_i(u)|\le C(1+|u|^2)$ due to (\Hh), we deduce from (\Pb) and (\La) in Lemma~2 that
$$\|u_i\|_{1,T} \le C(T)\Bigl[\|u_{i,0}\|_{C^1}+\|u_i\|_{0,T}^{1/2}\|f_i(u)\|_{0,T}^{1/2}\Bigr]
\le C(T)\bigl(1+\|u\|_{0,T}\bigr)^{3/2}.$$
Since $\nabla w_i= e^{-Kt}(1+\log(1+u_i))\nabla u_i$, it follows that
$$\|w_i\|_{1,T}\le (1+\log(1+\|u_i\|_{0,T}))\|u_i\|_{1,T}
\le C(T)\bigl(1+\|u\|_{0,T}\bigr)^{3/2}\log(2+\|u\|_{0,T}).$$
Combining this with (\Cc) and taking maximum over $i\in \{1,\dots,m\}$, we obtain
$$
\bigl(1+\|u\|_{0,T}\bigr)\log(1+\|u\|_{0,T})\le e^{KT}\|w\|_{0,T} 
\le C(T)\bigl(1+\|u\|_{0,T}\bigr)\Bigl(\log(2+\|u\|_{0,T})\Bigr)^{2/3},$$
hence $\|u\|_{0,T}\le C(T)$.
We conclude that $T_{max}=\infty$,
since $T_{max}<\infty$ would imply the blowup of $\|u\|_{0,T}$ as $T\to T_{max}$,
whereas $C(T)$ remains bounded for $T$ bounded.
\fin

\medskip

{\bf 3. Proof of Lemma~2.}
\medskip

In this proof, $C$ denotes a generic positive constant depending only on $\Omega, d, T$,
and remaining bounded for $T$ bounded.
\smallskip

(i) For $k>0$ to be chosen later, we note that $U$ solves
$U_t-d\Delta U+kU=g+kU$, hence 
$(\partial_t-d\Delta)(e^{kt}U)=e^{kt}(g+kU)$.
By the variation of constants formula, we deduce that
$$U(t)=e^{-kt}e^{td\Delta}U_0+\int_0^t e^{(t-s)d\Delta}e^{-k(t-s)}(g+kU)(s)\, ds,
\leqno(\Ld)$$
where $(e^{t\Delta})$ denotes the Neumann or the Cauchy heat semigroup.
We have the estimates 
$$\|e^{t\Delta}\psi\|_{C^1}\le C\|\psi\|_{C^1},\quad 0<t<T, \ \psi\in E_1,$$
and
$$\|e^{t\Delta}\psi\|_{C^1}\le Ct^{-1/2}\|\psi\|_{L^\infty},\quad 0<t<T, \ \psi\in BC(\overline\Omega),$$
(see e.g. [\Mo] and the references therein).
We deduce that
$$\|U(t)\|_{C^1} \le Ce^{-kt}\|U_0\|_{C^1}+C\int_0^t(t-s)^{-1/2}e^{-k(t-s)}\|(g+kU)(s)\|_\infty\, ds.$$
Using 
$$\int_0^\infty (t-s)^{-1/2}e^{-k(t-s)}\, ds=k^{-1/2}\int_0^\infty \tau^{-1/2}e^{-\tau}\, d\tau=Ck^{-1/2},$$
we obtain
$$\|U\|_{1,T} \le C \|U_0\|_{C^1}
+C\bigl(k^{-1/2}\|g\|_{0,T}+k^{1/2}\|U\|_{0,T}\bigr).$$
Inequality (\La) then follows by choosing $k=\|g\|_{0,T}\|U\|_{0,T}^{-1}$.
\smallskip

(ii) We have the estimates 
$$\|e^{t\Delta}\psi\|_{C^2}\le C\|\psi\|_{C^2},\quad 0<t<T, \ \psi\in E_2,$$
and
$$\|e^{t\Delta}\psi\|_{C^2}\le Ct^{-1/2}\|\psi\|_{C^1},\quad 0<t<T, \ \psi\in E_1$$
(see e.g. [\Mo] and the references therein).
It follows from (\Ld) that
$$\|U(t)\|_{C^2} \le Ce^{-kt}\|U_0\|_{C^2}+C\int_0^t(t-s)^{-1/2}e^{-k(t-s)}\|(g+kU)(s)\|_{C^1}\, ds.$$
By the argument in part (i), we deduce (\Lb).
Property (\Lc) then follows by combining (\La) and (\Lb).
\fin

\vskip 3mm
{\baselineskip=11pt \parindent=0.7cm

\font\rmn=cmr9
\font\sln=cmsl9
\font\rmb=cmbx8 scaled 1125 \rm

\rmn \eightpoint

\centerline{\bf REFERENCES}
\medskip\medskip

\item{[\BDS]}   M. Bisi, L. Desvillettes and G. Spiga,
Exponential convergence to equilibrium via Lyapounov functionals for reaction-diffusion equations 
arising from non reversible chemical kinetics,
{\sln ESAIM Math. Model. Numer. Anal.} 43 (2009), 151--172.

\smallskip

\item{[\CDF]} J.A. Ca\~nizo, L. Desvillettes and K. Fellner,
Improved duality estimates and applications to reaction-diffusion equations,
{\sln Comm. Partial Differential Equations} 39 (2014), 1185--1204.

\smallskip

\item{[\CGV]} M.C. Caputo, T. Goudon, A. Vasseur,
Solutions of the $4$-species quadratic reaction-diffusion system are bounded and $C^\infty$-smooth, in any space dimension,
Preprint arXiv:1709.05694 (2017).

\smallskip

\item{[\CV]} M.C. Caputo and A. Vasseur,
Global regularity of solutions to systems of reaction-diffusion with sub-quadratic growth in any dimension,
{\sln Comm. Partial Differential Equations} 34 (2009), 1228--1250.

\smallskip

\item{[\DF]} L. Desvillettes and K. Fellner,
Exponential decay toward equilibrium via entropy methods for reaction-diffusion equations
{\sln J. Math. Anal. Appl.} 319 (2006), 157--176.

\smallskip

\item{[\DFPV]} L. Desvillettes, K. Fellner, M. Pierre and J. Vovelle,
Global existence for quadratic systems of reaction-diffusion
{\sln Adv. Nonlinear Stud.} 7 (2007), 491--511.

\smallskip

\item{[\FLaa]}  K. Fellner and E.-H. Laamri,
Exponential decay towards equilibrium and global classical solutions for nonlinear reaction-diffusion systems,
{\sln J. Evol. Equ.} 16 (2016), 681--704.

\smallskip

\item{[\FPT]}  K. Fellner, W. Prager and B.Q. Tang
Exponential decay towards equilibrium and global classical solutions for nonlinear reaction-diffusion systems,
{\sln Kin\-et. Relat. Models} 10 (2017), 1055--1087.

\smallskip

\item{[\FT]}  K. Fellner and B.Q. Tang
Explicit exponential convergence to equilibrium for nonlinear reaction-diffusion systems with detailed balance condition
{\sln Nonlinear Anal.} 159 (2017), 145--180.

\smallskip

\item{[\GV]} T. Goudon and A. Vasseur,
Regularity analysis for systems of reaction-diffusion equations
{\sln Ann. Sci. \'Ec. Norm. Sup\'er.} (4) 43 (2010), 117--142.

\smallskip

\item{[\KaA]} J.I. Kanel,
The Cauchy problem for a system of semilinear parabolic equations with balance conditions,
{\sln Differentsial'nye Uravneniya} 20 (1984), 1753--1760
(English translation: {\sln Differential Equations} {\bf 20} (1984), 1260--1266).

\smallskip

\item{[\KaB]} J.I. Kanel,
Solvability in the large of a system of reaction-diffusion equations with the balance condition,
{\sln Differentsial'nye Uravneniya} 26 (1990), 448--458
(English translation: {\sln Differential Equations} {26} (1990), 331--339.

\smallskip

\item{[\MHM]} A. Mielke, J. Haskovec and P.A. Markowich,
On uniform decay of the entropy for reaction-diffusion systems,
{\sln J. Dynam. Differential Equations} 27 (2015), 897--928.

\smallskip

\item{[\Mo]} X. Mora, 
Semilinear parabolic equations define semiflows on $C^k$ spaces,
{\sln Trans. Amer. Math. Soc.} 278 (1983), 21--55.

\smallskip

\item{[\Pi]} M. Pierre,
Global existence in reaction-diffusion systems with control of mass: a survey,
{\sln Milan J. Math.} 78 (2010), 417--455.

\smallskip

\item{[\PSY]} M. Pierre, T. Suzuki and Y. Yamada,
Dissipative reaction diffusion systems with quadratic growth,
{\sln  Indiana Univ. Math. J.}  (2018), to appear (Preprint hal: 01671797).

\smallskip

\item{[\PSZ]} M. Pierre, T. Suzuki and R. Zou,
Asymptotic behavior of solutions to chemical reaction-diffusion systems,
{\sln J. Math. Anal. Appl.} 450 (2017), 152--168.

\smallskip

\item{[\QSb]} P. Quittner, Ph. Souplet,
Superlinear parabolic problems. Blow-up, global existence and steady states,
Birkhauser Advanced Texts, 2007, 584 p.+xi. 

\bye